\documentclass[reprint, twocolumn, aps]{revtex4-1}

\usepackage{natbib}
\usepackage{pdfx}
\usepackage{color}
\usepackage{amsfonts}
\usepackage{amsopn}
	
\usepackage{graphicx}
\usepackage{epstopdf}
\usepackage[caption=false]{subfig} 
\usepackage{algorithmic}
\usepackage{siunitx}
\usepackage{amsmath}
\usepackage{dsfont}
\usepackage{bm}
\usepackage{url}
\usepackage{algorithm,algorithmic,amssymb,letltxmacro}

\ifpdf
  \DeclareGraphicsExtensions{.eps,.pdf,.png,.jpg}
\else
  \DeclareGraphicsExtensions{.eps}
\fi
\newcommand{\xa}{\xi_{j}}
\newcommand{\ya}{\eta_{j}}
\newcommand{\xb}{\xi_{k}}
\newcommand{\yb}{\eta_{k}}
\newcommand{\xc}{\xi_{\ell}}
\newcommand{\yc}{\eta_{\ell}}
\newcommand{\PQ}{(\xb-\xa)^{2}+(\yb-\ya)^{2}}
\newcommand{\QR}{(\xc-\xb)^{2}+(\yc-\yb)^{2}}
\newcommand{\RP}{(\xa-\xc)^{2}+(\ya-\yc)^{2}}
\newcommand{\ed}[1]{\textcolor{black}{#1}}
\newcommand{\myname}{MC}
\begin{document}	

\title{A simple algorithm to find the L-curve corner \\ in the regularisation of ill-posed inverse problems}

\author{Alessandro Cultrera}\email{a.cultrera@inrim.it}
\affiliation{INRIM - Istituto Nazionale di Ricerca Metrologica, Strada delle Cacce, 91 --- 10135, Torino, Italy.}
\author{Luca Callegaro}
\affiliation{INRIM - Istituto Nazionale di Ricerca Metrologica, Strada delle Cacce, 91 --- 10135, Torino, Italy.}

\begin{abstract}
We propose a simple algorithm to locate the ``corner'' of an L-curve, a function often used to select the regularisation parameter for the solution of ill-posed inverse problems. The algorithm involves the Menger curvature of a circumcircle and the golden section search method. It efficiently finds the regularisation parameter value corresponding to the maximum positive curvature region of the L-curve. The algorithm is applied to some commonly available test problems and compared to the typical way of locating the l-curve corner by means of its analytical curvature. The application of the algorithm to the data processing of an electrical resistance tomography experiment on thin conductive films is also reported. 

\vspace{5mm}

Now published at \emph{IOP SciNotes} 1 (2020) 025004,~ \url{https://doi.org/10.1088/2633-1357/abad0d}
\end{abstract}

\maketitle

\section{Introduction}	
The solution $\bm{\hat x}$ of an ill-posed inverse problem is often searched by means of a regularized least squares functional of the type 
\begin{equation}
\begin{aligned}
\label{eqn:funct}
\bm{\hat x}_{\lambda}= \textup{arg }  \underset{\bm{x}}{\textup{min}} \left \{ ||\bm{A}\bm{x}-\bm{b}||^{2}+\lambda \bm{R}(\bm{x})\right \},\quad \lambda \in \mathds{R}, \quad \lambda \geq 0 
\end{aligned}
\end{equation}
 %
where $\bm{Ax}-\bm{b}$ is the vector of residuals between the experimental data vector $\bm{b}$ and the reconstructed data $\bm{A}\bm{x}$ for a given $\bm{x}$. The regularisation term $\bm{R}(\bm{x})$ renders the problem less sensitive to the noise of $\bm{b}$ and find a stable solution. 
$\bm{R}(\bm{x})$ represents a cost function, which usually includes prior information about the solution. The scalar factor $\lambda$ is the \emph{regularisation parameter}, is a weighing factor of $\bm{R}(\bm{x})$. The choice of $\lambda$ is crucial for a meaningful solution.
As an example, we consider the regularisation method of Tikhonov~\cite{tikhonov2013numerical}, in which $\bm{R}(\bm{x}) =  ||\bm{x}||^{2}$.
%
%
Several methods~(see \cite[\textsection.~7]{hansen1998rank}) have been developed in order to find an optimal tuning of $\lambda$ for a given problem. Of particular interest is the L-curve method \cite[\textsection.~7.5]{hansen1998rank} \cite{hansen1992analysis}, which is one of the best-known heuristic methods for the selection of $\lambda$. 
The L-curve is two-dimensional, parametric in $\lambda$, defined by points with cartesian coordinates
\begin{equation}
\label{eqn:lcurve}
P(\lambda)=(\xi(\lambda), \eta(\lambda))\rightarrow\begin{cases}
& \xi(\lambda)= \log||\bm{Ax}-\bm{b}||^2\\
&\eta(\lambda)= \log||\bm{x}||^2
\end{cases}
\end{equation}
The point of maximum positive curvature $P(\lambda_\mathrm{opt})$, the ``corner'', can be associated to the optimal reconstruction parameter, say $\lambda_\mathrm{opt}$.  The underlying concept is that the ``corner'' represents a compromise between the fitting to the data and the amount of regularisation applied to the problem~\cite{hansen1993use}.
Numerical search algorithms have been proposed for the estimation of $\lambda_\mathrm{opt}$; among them, we mention the splines method~\cite{hansen1992analysis,hansen2007adaptive}, the triangle method~\cite{castellanos2002triangle} and the L-ribbon method~\cite{calvetti1999estimation}.
The adoption of the L-curve approach to deal with diverse ill-posed inverse problems is an ongoing research topic~\cite{choi2019interpretation, xu2016extended}. 
Here we propose an alternative method, and its very simple implementation, to locate the L-curve corner. It is based on an iterative estimation of the local curvature of the L-curve from three sampled points with an update rule based on the golden section search. The method has a small computational effort since it reduces the number points of the L-curve explicitly computed. The following gives a description of the algorithm and its application on both typical test problems, and a reconstruction problem of electrical resistance tomography.
%
\section{Algorithm}
\label{sec:alg}
The algorithm~\ref{alg} is written in pseudo-code. Algorithm~\ref{alg} calls two functions. Function P = \texttt{l\_curve\_P($\lambda$)} is based on the the specific regularisation problem being solved; it is assumed that at each call, given as input the regularisation parameter $\lambda$ it solves the system (1) and provides as output the point $P(\lambda)$, i.e. the coordinates $\xi(\lambda)$ and $\eta(\lambda)$ of the L-curve. The function $C_k = \texttt{menger($P_j$, $P_k$, $P_\ell$)}$ is defined below in Sec.~\ref{sec:curv}.
The algorithm is iterative and identifies the estimate $\lambda_\mathrm{opt}$, in the following $\lambda_\mathrm{\myname}$, by means of the definition of curvature given in Sec.~\ref{sec:curv} and the golden section search method, described in Sec.~\ref{sec:gss}. ``MC'' stays for ``Menger Curvature''.
\subsection{Curvature}
\label{sec:curv}
The function $C_k = \texttt{menger($P_j$, $P_k$, $P_\ell$)}$ is based on the definition of the curvature of a circle by three points given by Menger~\cite{menger1930untersuchungen,pajot2002analytic}. In our case three values $\lambda_{j}<\lambda_{k}<\lambda_{\ell}$ of the regularisation parameter identify three points $P(\lambda_{j})$, $P(\lambda_{k})$ and $P(\lambda_{\ell})$ on the L-curve. 
We follow the notation of~\eqref{eqn:lcurve} for the coordinates of a generic point $P(\lambda)$. For notational simplicity we make the substitution:
\begin{equation}
\begin{aligned}
\label{eqn:subst}
\xi(\lambda_i)&\rightarrow \xi_i,\\
\eta(\lambda_i)&\rightarrow \eta_i,\\
P(\lambda_i)&\rightarrow P_i.
\end{aligned}
\end{equation}
We define a signed curvature $C_{k}$ of the circumcircle as
\begin{equation}
\label{eqn:menga}
\begin{aligned}
C_k = \frac{2\cdot(\xa\yb+\xb\yc+\xc\ya-\xa\yc-\xb\ya-\xc\yb)}
{\Big(\overline{P_{j}P_{k}}\cdot \overline{P_{k}P_{\ell}}\cdot \overline{P_{\ell}P_{j})}\Big)^{1/2}}, 
\end{aligned}
\end{equation}
where
\begin{equation}
\begin{aligned}
\label{eqn:defs}
\overline{P_{j}P_{k}}&=\PQ,\\
\overline{P_{k}P_{\ell}}&=\QR,\\
\overline{P_{\ell}P_{j}}&=\RP, 
\end{aligned}
\end{equation}
are the euclidean distances between the sampled L-curve points. Note that we choose to index the curvature with the intermediate index ($k$) of the three points.
\subsection{Golden Section Search}
\label{sec:gss}
The algorithm is initialized by assigning the search interval $[\lambda_1,\lambda_4]$. Two other values $\lambda_{2}$ and $\lambda_{3}$ are calculated following the golden section method; the calculation is done on the exponents of $\lambda$ (given $\lambda_i = 10^{x_i}$) to maintain a uniform spacing along the many orders of magnitude covered,
\begin{equation}
\begin{aligned}
\label{eqn:goldenratio}
x_{2} &= (x_{4}+\varphi\cdot x_{1})/(1+\varphi), \\
x_{3} &= x_{1}+(x_{4}-x_{2}),\\
\end{aligned}
\end{equation}
where $\varphi = (1+\sqrt{5})/2$ is the golden section~\cite{ kiefer1953sequential}.  %
Four values of $\lambda$ define four points on the L-curve and allow to calculate two curvatures, $C_{2}$ from $\{P(\lambda_{1})$, $P(\lambda_{2})$, $P(\lambda_{3})\}$ and $C_{3}$ from $\{P(\lambda_{2})$, $P(\lambda_{3})$, $P(\lambda_{4})\}$. 
The curvatures $C_{2}$ and $C_{3}$ are compared; consistent reassignment and recalculation are done in order to work at each iteration with four points $P(\lambda_1)\ldots P(\lambda_4)$.  
The algorithm terminates when the search interval 
$[\lambda_1,\lambda_4]$ is smaller than a specified threshold $\epsilon$ and returns $\lambda_\mathrm{\myname}$.
 
It may happen that the curvature $C_3$ associated to the right-hand circle is negative at the initial stage of the search, since $C_k$ is defined with sign in~\eqref{eqn:menga}. By definition, the corner corresponds to a positive curvature and it lays on the left-side of the plot. Hence, the algorithm performs a check, and while $C_{3} < 0$ the search extreme $\lambda_{1}$ is kept fixed, $\lambda_{4}$ is shifted toward smaller values and   $\lambda_{2}$ and  $\lambda_{3}$ are recalculated. The condition on $C_3$ is strong enough that even in case of both negative curvatures it guarantees the convergence towards the corner.

Some considerations: (a) according to the golden section search method, the algorithm needs to recalculate only one $P(\lambda)$ at each iteration (except for the first iteration), the other can be simply reassigned; this limits the calculation effort; 
(b) as $P(\lambda_{1})$ and $P(\lambda_{4})$ are distant at the first iterations, $C_2$ and $C_3$ are just rough approximations of the curvature of the L-curve in different regions, but become more accurate as the distance between the search extremes decreases.
\begin{algorithm}[H]
\caption{L-curve corner search}
\label{alg}
\begin{algorithmic}[1]
\STATE{Initialize  $\lambda_{1}$ and $\lambda_{4}$;}
	\COMMENT{search extremes}
\STATE{Assign $\epsilon$;}
	\COMMENT{termination threshold}
\STATE{$\varphi \leftarrow (1+\sqrt{5})/2$;}
	\COMMENT{golden section}
\STATE{$\lambda_{2} \leftarrow \ed{10^{(x_{4}+\varphi\cdot x_{1})/(1+\varphi)}}$;}
\STATE{$\lambda_{3} \leftarrow \ed{10^{x_{1}+(x_{4}-x_{2})}}$;}
 \FOR{$i=1$ to $4$}
 	\STATE{$P_{i}\leftarrow$\texttt{l\_curve\_P}($\lambda_i$);}
		\COMMENT{\texttt{l\_curve\_P} returns~\eqref{eqn:lcurve}}
 \ENDFOR
\REPEAT
\STATE{$C_2 \leftarrow$\texttt{menger}($P_1$,$P_2$,$P_3$);}
	\COMMENT{\texttt{menger} calls~\eqref{eqn:menga}}
\STATE{$C_3 \leftarrow$\texttt{menger}($P_2$,$P_3$,$P_4$);}
	
	\REPEAT 
		\STATE{$\lambda_{4} \leftarrow \lambda_{3}$;\quad $P_{4} \leftarrow P_{3}$;}
		\STATE{$\lambda_{3} \leftarrow \lambda_{2}$;\quad $P_{3} \leftarrow P_{2}$;}
		\STATE{$\lambda_{2} \leftarrow \ed{10^{(x_{4}+\varphi\cdot x_{1})/(1+\varphi)}}$;}
		\STATE{$P_{2}\leftarrow$\texttt{l\_curve\_P}($\lambda_2$);}
		\STATE{$C_3 \leftarrow$\texttt{menger}($P_2$,$P_3$,$P_4$);}
	\UNTIL{$C_3 >0$}

\IF{$C_{2} > C_{3}$}

	\STATE{$\lambda\leftarrow\lambda_{2}$;}
		\COMMENT{store $\lambda$}
	\STATE{$\lambda_{4} \leftarrow \lambda_{3}$; \quad $P_{4} \leftarrow P_{3}$;}
	\STATE{$\lambda_{3} \leftarrow \lambda_{2}$; \quad $P_{3} \leftarrow P_{2}$;}
	\STATE{$\lambda_{2} \leftarrow \ed{10^{(x_{4}+\varphi\cdot x_{1})/(1+\varphi)}}$;}
	\STATE{$P_{2}\leftarrow$\texttt{l\_curve\_P}($\lambda_2$);}	
		\COMMENT{only $P_2$ is recalculated}
\ELSE
	\STATE{$\lambda\leftarrow\lambda_{3}$}
	\STATE{$\lambda_{1} \leftarrow \lambda_{2}$; \quad $P_{1} \leftarrow P_{2}$;}
	\STATE{$\lambda_{2} \leftarrow \lambda_{3}$; \quad $P_{2} \leftarrow P_{3}$;}
	\STATE{$\lambda_{3} \leftarrow \ed{10^{x_{1}+(x_{4}-x_{2})}}$;}
	\STATE{$P_{3}\leftarrow$\texttt{l\_curve\_P}($\lambda_3$);}
		\COMMENT{only $P_3$ is recalculated}	
\ENDIF
\UNTIL{ $(\lambda_{4}-\lambda_{1})/\lambda_{4}<\epsilon$}
\RETURN $\lambda_\mathrm{\myname}\leftarrow\lambda$
\end{algorithmic} 
\end{algorithm}
\begin{figure*}
\centering 
\subfloat[Iteration 1]{\label{fig:lcurves1}\includegraphics[width=4cm]{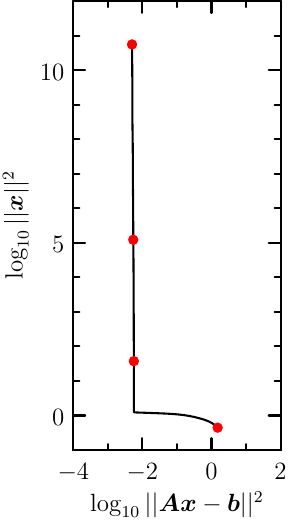}} 
\subfloat[Iteration 2]{\label{fig:lcurves2}\includegraphics[width=4cm]{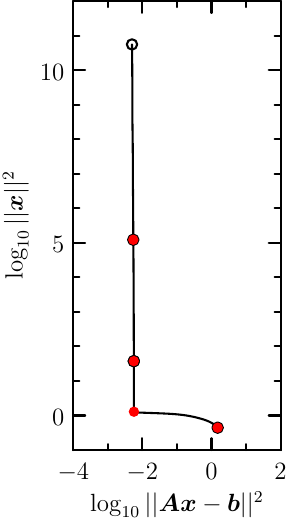}} 
\subfloat[Iteration 3]{\label{fig:lcurves3}\includegraphics[width=4cm]{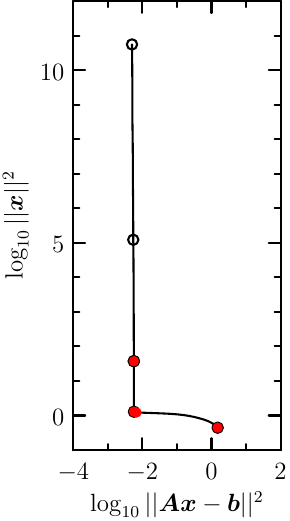}} 
\caption{The algorithm at \ed{the first three} iterations. The reference L-curve is reported as a solid line. Solid circles represent the points $P(\lambda)$ being evaluated at the labeled iteration. Empty circles represent the points evaluated at past iterations. }
\label{fig:lcurves}
\end{figure*}
\section{Application to test problems}
\label{sec:test}
\ed{We tested the algorithm on small demonstrative problems, some (\texttt{baart}, \texttt{blur}, \texttt{shaw} and \texttt{spike})  chosen from the function library \texttt{Regularisation Tools} (\texttt{RT}), implemented in \texttt{MATLAB}~\cite{hansen1994regularization}. This library is also employed to implement a function of the algorithm (\texttt{L\_curve\_P($\lambda$)}) which evaluates a single point of the L-curve for a given $\lambda$. Algorithm~\ref{alg} is implemented in \texttt{MATLAB} as well.}

\begin{figure}
\centering 
\includegraphics[width=0.9\columnwidth]{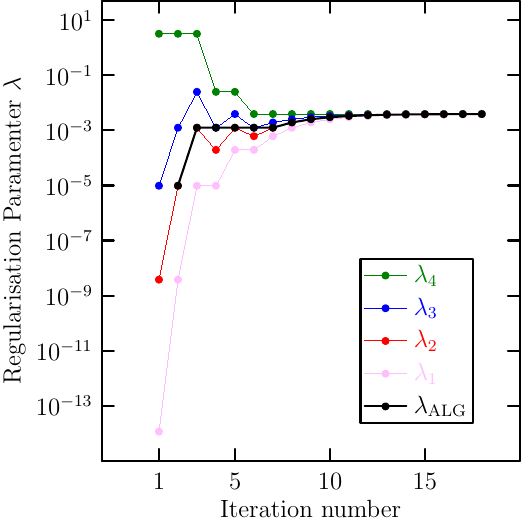} 
\caption{Behaviour of the algorithm versus iteration number for the \texttt{baart} problem with a threshold $\epsilon = 1\%$. The last point (iteration 18) corresponds to $\lambda_\mathrm{\myname}$.}
\label{fig:conv}
\end{figure}

\ed{The application of algorithm~\ref{alg} to the problem \texttt{baart} is shown explicitly in the following. This problem represents the discretization of a Fredholm integral equation of first kind of order $n$. The matrix $\bm{A}$ in \eqref{eqn:funct} is therefore $n \times n$. The chosen size of the problem is $n = 32$. In this example we added random noise of relative standard deviation of $10^{-3}$ to the exact data. The corner of the L-curve generated by this problem is located with both  algorithm~\ref{alg} and the \texttt{L\_corner} routine from \texttt{RT}~\footnote{We tested our algorithm with noise relative standard deviation levels over a wide range, from $10^{-10}$ to $10^{-1}$. Our results always matched with negligible deviation the algorithm of \texttt{Regularisation Tools}, taken as reference.}.
Figure~\ref{fig:lcurves}  shows the first three iterations of the algorithm, and displays also a full L-curve obtained by dense sampling of \texttt{L\_curve\_P($\lambda$)} as a reference. Empty circles represent points visited at previous iterations, while filled circles represent the four points $P_1 \ldots P_4$ of the given iteration.
The algorithm runs by choosing as initial search extremes the default choice of the \texttt{L\_corner} routine ($\lambda_{1} = 10^{-14}$ and $\lambda_{4} = 10^{-1}$).
%
%
Running the algorithm on the other three mentioned problems gives similar results in term of accuracy compared to the native \texttt{L\_corner} routine of the \texttt{RT} library. The optimal regularisation parameter obtained with this routine is called $\lambda_\mathrm{RT}$ in the following.
Table~\ref{tab:test} summarizes the results of solving the four test problems with algorithm~\ref{alg} and with \texttt{RT}'s function \texttt{L\_corner}. $\lambda_\mathrm{\myname}$ is the optimal regularisation parameter returned by our algorithm while $\lambda_\mathrm{RT}$ is the one returned by the \texttt{RT} routine. The \texttt{MATLAB} profiler was used to get the corresponding net timing~\footnote{Eventual plotting time not considered.} of algorithm~\ref{alg} ($t_\mathrm{\myname}$), and the \texttt{L\_corner} routine ($t_\mathrm{RT}$).} 
%
%
%
\begin{table*}
\renewcommand{\arraystretch}{1.3}
\caption{\ed{Comparison between algorithm~\ref{alg} and the analytic curvature approach on test problems.}}
\label{tab:test}
\centering
\begin{tabular}{lccccc}
\hline 
\bfseries Problem & $\lambda_\mathrm{\myname}$  &  $\lambda_\mathrm{RT}$ &  $t_\mathrm{\myname}$  (\si{\milli\second}) & $t_\mathrm{RT}$ (\si{\milli\second}) & iterations  \\
\hline
\texttt{baart(32)} 	& $3.92\times10^{-3}$	& $4.02\times10^{-3}$	& \num{73}		&\num{469} 	&\num{18} \\
\texttt{blur(16,4,5)} 	& $3.18\times10^{-4}$	& $3.20\times10^{-4}$	& \num{82}		&\num{459}	&\num{15}\\
\texttt{shaw(32)} 	& $8.65\times10^{-4}$	& $8.28\times10^{-4}$	& \num{69}		&\num{473}	&\num{17}\\
\texttt{spike(32,5)} 	& $1.65\times10^{-4}$	& $1.60\times10^{-4}$	& \num{61}		&\num{461}	&\num{17}\\
\hline

\end{tabular}
\end{table*}
%
Fig.~\ref{fig:conv} shows the evolution of the algorithm towards convergence. 

As a side note, a similar implementation of the presented algorithm could be made also using a Fibonacci search to pick the $x_i$ in \eqref{eqn:goldenratio}. In fact the Fibonacci search interval reduction ratio converges to the golden section very quickly~\cite{ kiefer1953sequential}.
\section{Application to Electrical Resistance Tomography}
\label{sec:experiments}
\begin{figure*}
\centering 
\subfloat[Scheme of the sample and contacts. Colors represent nominal conductivity. The sample was patterned to remove the conductive coating and obtain non-conductive features (black)]{\label{fig:reca}\includegraphics[width=0.8\columnwidth]{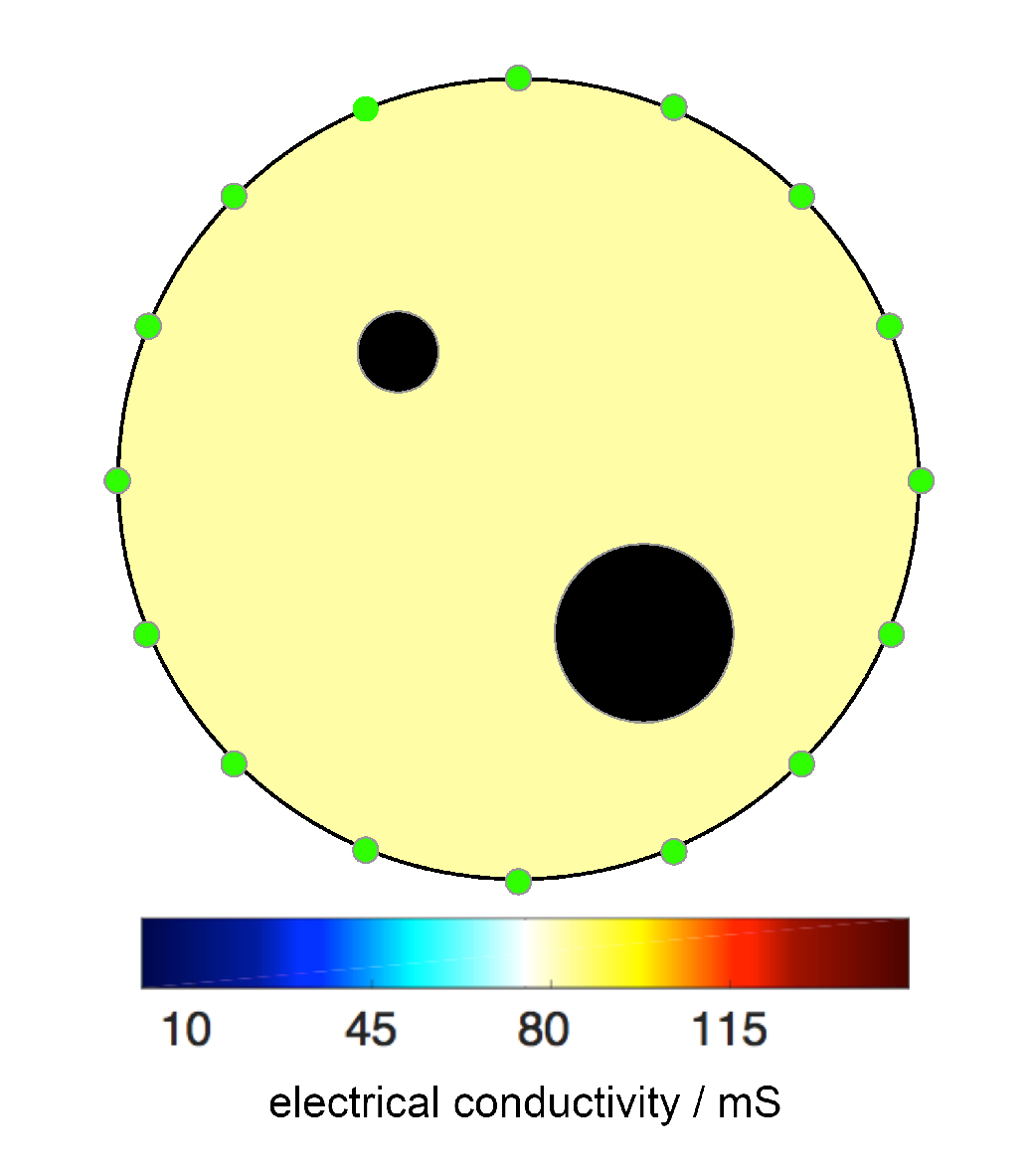}} 
\qquad
\subfloat[Conductivity map reconstructed from experimental data with $\lambda=1.764\cdot10^{-6}$.]{\label{fig:recb}\includegraphics[width=0.8\columnwidth]{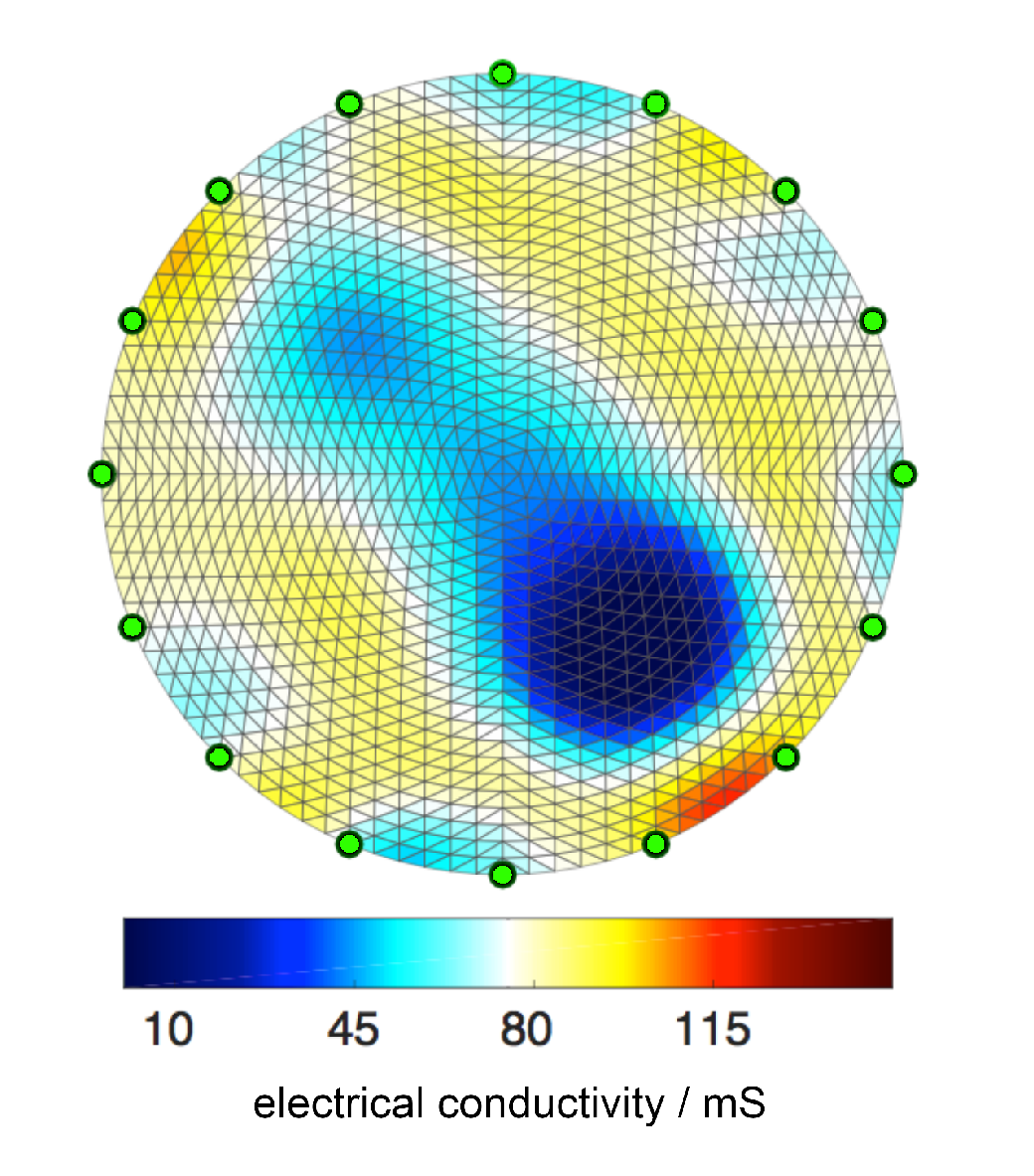}} 
\caption{Schematic of the conductive sample (a) and ERT conductivity map obtained with $\lambda_\mathrm{\myname}$ from our algorithm (b).}
\label{fig:rec}
\end{figure*}
\begin{figure}
\centering 
\subfloat[Iteration 17]{\label{fig:overviewERT}\includegraphics[width=0.33\columnwidth]{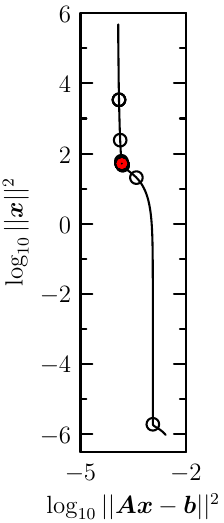}} 
\subfloat[``Corner'' detail at iteration 17]{\label{fig:ERTzoom}\includegraphics[width=0.66\columnwidth]{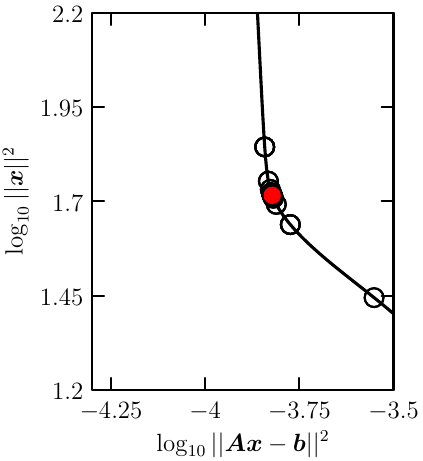}} 
\caption{L-curve obtained from experimental data and the points corresponding to the last iteration of the algorithm (a). Detail of the corner (b).}
\label{fig:ert}
\end{figure}

The following shows the application of algorithm~\ref{alg} to electrical resistance tomography (ERT)~\cite{seo2013electrical}. 
In this experiment we used a patterned tin-oxide conductive sample of circular geometry, with electrical contacts on its boundary (see Fig.~\ref{fig:reca}). 
Four-terminal resistance measurements are performed with a scanning setup; the measurements are the elements of the data vector $\bm{b}$. A detailed description of the experiment is given in~\cite{cultrera2019mapping, cultrera2016electrical}.  
The ERT problem solution is obtained by solving a discretized Laplace equation with Tikhonov regularisation, a formulation compatible with the calculation of a continuous L-curve. 
\texttt{EIDORS}~\cite{adler2005eidors} routines are used to generate a two-dimensional circular mesh (2304 elements) with 16 contact points at the boundary (corresponding to a $\bm{b}$ of size 208), to discretize the Laplace equation and obtain matrix $\bm{A}$. The reconstructed image shown in Fig.~\ref{fig:recb}.
Figure~\ref{fig:ert} reports the main results of the application of algorithm~\ref{alg} to ERT experimental data. Figure~\ref{fig:overviewERT} shows the L-curve and the last iteration of algorithm~\ref{alg}; figure~\ref{fig:ERTzoom} the detail of the corner. 
The optimal regularisation parameter returned by the algorithm with $\epsilon = 1\%$ is $\lambda_\mathrm{\myname} = 1.8\cdot10^{-6}$. The relative difference between $\lambda_\mathrm{\myname}$ and $\lambda_\mathrm{RT}$ is negligible (1 part in $10^{13}$). 
\section{Conclusions}
\label{sec:conclusions}
The proposed algorithm allows, given an inverse problem having the form \eqref{eqn:funct}, the determination of the regularisation parameter $\lambda_\mathrm{\myname}$ corresponding to the maximum positive curvature of the L-curve. The algorithm is designed for maximum simplicity of implementation on already existing solvers. On both test problems or in a  real electrical resistance tomography problem, convergence is achieved in less than 20 iterations. Compared to a common routine for the location of the L-curve corner such as \texttt{Regularisation Tools} the present algorithm returns strongly compatible results with a reduced calculation effort.
\section*{Acknowledgements}
The work has been developed within the Joint Research Project 16NRM01 GRACE: Developing electrical characterisation methods for future graphene electronics. This project has received funding from the EMPIR programme co-financed by the Participating States and from the European Union’s Horizon 2020 research and innovation programme.

\newpage

\begin{thebibliography}{19}%
	\makeatletter
	\providecommand \@ifxundefined [1]{%
		\@ifx{#1\undefined}
	}%
	\providecommand \@ifnum [1]{%
		\ifnum #1\expandafter \@firstoftwo
		\else \expandafter \@secondoftwo
		\fi
	}%
	\providecommand \@ifx [1]{%
		\ifx #1\expandafter \@firstoftwo
		\else \expandafter \@secondoftwo
		\fi
	}%
	\providecommand \natexlab [1]{#1}%
	\providecommand \enquote  [1]{``#1''}%
	\providecommand \bibnamefont  [1]{#1}%
	\providecommand \bibfnamefont [1]{#1}%
	\providecommand \citenamefont [1]{#1}%
	\providecommand \href@noop [0]{\@secondoftwo}%
	\providecommand \href [0]{\begingroup \@sanitize@url \@href}%
	\providecommand \@href[1]{\@@startlink{#1}\@@href}%
	\providecommand \@@href[1]{\endgroup#1\@@endlink}%
	\providecommand \@sanitize@url [0]{\catcode `\\12\catcode `\$12\catcode
		`\&12\catcode `\#12\catcode `\^12\catcode `\_12\catcode `\%12\relax}%
	\providecommand \@@startlink[1]{}%
	\providecommand \@@endlink[0]{}%
	\providecommand \url  [0]{\begingroup\@sanitize@url \@url }%
	\providecommand \@url [1]{\endgroup\@href {#1}{\urlprefix }}%
	\providecommand \urlprefix  [0]{URL }%
	\providecommand \Eprint [0]{\href }%
	\providecommand \doibase [0]{http://dx.doi.org/}%
	\providecommand \selectlanguage [0]{\@gobble}%
	\providecommand \bibinfo  [0]{\@secondoftwo}%
	\providecommand \bibfield  [0]{\@secondoftwo}%
	\providecommand \translation [1]{[#1]}%
	\providecommand \BibitemOpen [0]{}%
	\providecommand \bibitemStop [0]{}%
	\providecommand \bibitemNoStop [0]{.\EOS\space}%
	\providecommand \EOS [0]{\spacefactor3000\relax}%
	\providecommand \BibitemShut  [1]{\csname bibitem#1\endcsname}%
	\let\auto@bib@innerbib\@empty
	\bibitem [{\citenamefont {Tikhonov}\ \emph {et~al.}(2013)\citenamefont
		{Tikhonov}, \citenamefont {Goncharsky}, \citenamefont {Stepanov},\ and\
		\citenamefont {Yagola}}]{tikhonov2013numerical}%
	\BibitemOpen
	\bibfield  {author} {\bibinfo {author} {\bibfnamefont {A.~N.}\ \bibnamefont
			{Tikhonov}}, \bibinfo {author} {\bibfnamefont {A.}~\bibnamefont
			{Goncharsky}}, \bibinfo {author} {\bibfnamefont {V.}~\bibnamefont
			{Stepanov}}, \ and\ \bibinfo {author} {\bibfnamefont {A.~G.}\ \bibnamefont
			{Yagola}},\ }\href@noop {} {\emph {\bibinfo {title} {Numerical methods for
				the solution of ill-posed problems}}}\ (\bibinfo  {publisher} {Springer
		Science \& Business Media, New York, US},\ \bibinfo {year}
	{2013})\BibitemShut {NoStop}%
	\bibitem [{\citenamefont {Hansen}(1998)}]{hansen1998rank}%
	\BibitemOpen
	\bibfield  {author} {\bibinfo {author} {\bibfnamefont {P.~C.}\ \bibnamefont
			{Hansen}},\ }\href@noop {} {\emph {\bibinfo {title} {Rank-deficient and
				discrete ill-posed problems: numerical aspects of linear inversion}}}\
	(\bibinfo  {publisher} {Siam},\ \bibinfo {year} {1998})\BibitemShut {NoStop}%
	\bibitem [{\citenamefont {Hansen}(1992)}]{hansen1992analysis}%
	\BibitemOpen
	\bibfield  {author} {\bibinfo {author} {\bibfnamefont {P.~C.}\ \bibnamefont
			{Hansen}},\ }\href@noop {} {\bibfield  {journal} {\bibinfo  {journal} {SIAM
				Rev.}\ }\textbf {\bibinfo {volume} {34}},\ \bibinfo {pages} {561} (\bibinfo
		{year} {1992})}\BibitemShut {NoStop}%
	\bibitem [{\citenamefont {Hansen}\ and\ \citenamefont
		{O'Leary}(1993)}]{hansen1993use}%
	\BibitemOpen
	\bibfield  {author} {\bibinfo {author} {\bibfnamefont {P.~C.}\ \bibnamefont
			{Hansen}}\ and\ \bibinfo {author} {\bibfnamefont {D.~P.}\ \bibnamefont
			{O'Leary}},\ }\href@noop {} {\bibfield  {journal} {\bibinfo  {journal} {SIAM
				J. Sci. Comput.}\ }\textbf {\bibinfo {volume} {14}},\ \bibinfo {pages} {1487}
		(\bibinfo {year} {1993})}\BibitemShut {NoStop}%
	\bibitem [{\citenamefont {Hansen}\ \emph {et~al.}(2007)\citenamefont {Hansen},
		\citenamefont {Jensen},\ and\ \citenamefont
		{Rodriguez}}]{hansen2007adaptive}%
	\BibitemOpen
	\bibfield  {author} {\bibinfo {author} {\bibfnamefont {P.~C.}\ \bibnamefont
			{Hansen}}, \bibinfo {author} {\bibfnamefont {T.~K.}\ \bibnamefont {Jensen}},
		\ and\ \bibinfo {author} {\bibfnamefont {G.}~\bibnamefont {Rodriguez}},\
	}\href@noop {} {\bibfield  {journal} {\bibinfo  {journal} {J. Comput. Appl.
				Math.}\ }\textbf {\bibinfo {volume} {198}},\ \bibinfo {pages} {483} (\bibinfo
		{year} {2007})}\BibitemShut {NoStop}%
	\bibitem [{\citenamefont {Castellanos}\ \emph {et~al.}(2002)\citenamefont
		{Castellanos}, \citenamefont {G{\'o}mez},\ and\ \citenamefont
		{Guerra}}]{castellanos2002triangle}%
	\BibitemOpen
	\bibfield  {author} {\bibinfo {author} {\bibfnamefont {J.~L.}\ \bibnamefont
			{Castellanos}}, \bibinfo {author} {\bibfnamefont {S.}~\bibnamefont
			{G{\'o}mez}}, \ and\ \bibinfo {author} {\bibfnamefont {V.}~\bibnamefont
			{Guerra}},\ }\href@noop {} {\bibfield  {journal} {\bibinfo  {journal} {Appl.
				Numer. Math.}\ }\textbf {\bibinfo {volume} {43}},\ \bibinfo {pages} {359}
		(\bibinfo {year} {2002})}\BibitemShut {NoStop}%
	\bibitem [{\citenamefont {Calvetti}\ \emph {et~al.}(1999)\citenamefont
		{Calvetti}, \citenamefont {Golub},\ and\ \citenamefont
		{Reichel}}]{calvetti1999estimation}%
	\BibitemOpen
	\bibfield  {author} {\bibinfo {author} {\bibfnamefont {D.}~\bibnamefont
			{Calvetti}}, \bibinfo {author} {\bibfnamefont {G.~H.}\ \bibnamefont {Golub}},
		\ and\ \bibinfo {author} {\bibfnamefont {L.}~\bibnamefont {Reichel}},\
	}\href@noop {} {\bibfield  {journal} {\bibinfo  {journal} {BIT Numer. Math.}\
		}\textbf {\bibinfo {volume} {39}},\ \bibinfo {pages} {603} (\bibinfo {year}
		{1999})}\BibitemShut {NoStop}%
	\bibitem [{\citenamefont {Choi}\ \emph {et~al.}(2019)\citenamefont {Choi},
		\citenamefont {Shin}, \citenamefont {Ji}, \citenamefont {Kim}, \citenamefont
		{Son}, \citenamefont {Lee}, \citenamefont {Kim}, \citenamefont {Lee},\ and\
		\citenamefont {Yoon}}]{choi2019interpretation}%
	\BibitemOpen
	\bibfield  {author} {\bibinfo {author} {\bibfnamefont {M.-B.}\ \bibnamefont
			{Choi}}, \bibinfo {author} {\bibfnamefont {J.}~\bibnamefont {Shin}}, \bibinfo
		{author} {\bibfnamefont {H.-I.}\ \bibnamefont {Ji}}, \bibinfo {author}
		{\bibfnamefont {H.}~\bibnamefont {Kim}}, \bibinfo {author} {\bibfnamefont
			{J.-W.}\ \bibnamefont {Son}}, \bibinfo {author} {\bibfnamefont {J.-H.}\
			\bibnamefont {Lee}}, \bibinfo {author} {\bibfnamefont {B.-K.}\ \bibnamefont
			{Kim}}, \bibinfo {author} {\bibfnamefont {H.-W.}\ \bibnamefont {Lee}}, \ and\
		\bibinfo {author} {\bibfnamefont {K.~J.}\ \bibnamefont {Yoon}},\ }\href
	{\doibase 10.1007/s11837-019-03762-8} {\bibfield  {journal} {\bibinfo
			{journal} {JOM}\ }\textbf {\bibinfo {volume} {71}},\ \bibinfo {pages} {3825}
		(\bibinfo {year} {2019})}\BibitemShut {NoStop}%
	\bibitem [{\citenamefont {Xu}\ \emph {et~al.}(2016)\citenamefont {Xu},
		\citenamefont {Pei},\ and\ \citenamefont {Dong}}]{xu2016extended}%
	\BibitemOpen
	\bibfield  {author} {\bibinfo {author} {\bibfnamefont {Y.}~\bibnamefont
			{Xu}}, \bibinfo {author} {\bibfnamefont {Y.}~\bibnamefont {Pei}}, \ and\
		\bibinfo {author} {\bibfnamefont {F.}~\bibnamefont {Dong}},\ }\href {\doibase
		10.1088/0957-0233/27/11/114002} {\bibfield  {journal} {\bibinfo  {journal}
			{Meas. Sci. Technol.}\ }\textbf {\bibinfo {volume} {27}},\ \bibinfo {pages}
		{114002} (\bibinfo {year} {2016})}\BibitemShut {NoStop}%
	\bibitem [{\citenamefont {Menger}(1930)}]{menger1930untersuchungen}%
	\BibitemOpen
	\bibfield  {author} {\bibinfo {author} {\bibfnamefont {K.}~\bibnamefont
			{Menger}},\ }\href@noop {} {\bibfield  {journal} {\bibinfo  {journal} {Math.
				Ann.}\ }\textbf {\bibinfo {volume} {103}},\ \bibinfo {pages} {466} (\bibinfo
		{year} {1930})}\BibitemShut {NoStop}%
	\bibitem [{\citenamefont {Pajot}(2002)}]{pajot2002analytic}%
	\BibitemOpen
	\bibfield  {author} {\bibinfo {author} {\bibfnamefont {H.}~\bibnamefont
			{Pajot}},\ }\href@noop {} {\emph {\bibinfo {title} {Analytic capacity,
				rectifiability, Menger curvature and Cauchy integral}}}\ (\bibinfo
	{publisher} {Springer Science \& Business Media},\ \bibinfo {year}
	{2002})\BibitemShut {NoStop}%
	\bibitem [{\citenamefont {Kiefer}(1953)}]{kiefer1953sequential}%
	\BibitemOpen
	\bibfield  {author} {\bibinfo {author} {\bibfnamefont {J.}~\bibnamefont
			{Kiefer}},\ }\href@noop {} {\bibfield  {journal} {\bibinfo  {journal} {P. Am.
				Math. Soc.}\ }\textbf {\bibinfo {volume} {4}},\ \bibinfo {pages} {502}
		(\bibinfo {year} {1953})}\BibitemShut {NoStop}%
	\bibitem [{\citenamefont {Hansen}(1994)}]{hansen1994regularization}%
	\BibitemOpen
	\bibfield  {author} {\bibinfo {author} {\bibfnamefont {P.~C.}\ \bibnamefont
			{Hansen}},\ }\href@noop {} {\bibfield  {journal} {\bibinfo  {journal} {Numer.
				Algorithms}\ }\textbf {\bibinfo {volume} {6}},\ \bibinfo {pages} {1}
		(\bibinfo {year} {1994})}\BibitemShut {NoStop}%
	\bibitem [{Note1()}]{Note1}%
	\BibitemOpen
	\bibinfo {note} {We tested our algorithm with noise relative standard
		deviation levels over a wide range, from $10^{-10}$ to $10^{-1}$. Our results
		always matched with negligible deviation the algorithm of \protect \texttt
		{Regularisation Tools}, taken as reference.}\BibitemShut {Stop}%
	\bibitem [{Note2()}]{Note2}%
	\BibitemOpen
	\bibinfo {note} {Eventual plotting time not considered.}\BibitemShut {Stop}%
	\bibitem [{\citenamefont {Seo}\ and\ \citenamefont
		{Woo}(2013)}]{seo2013electrical}%
	\BibitemOpen
	\bibfield  {author} {\bibinfo {author} {\bibfnamefont {J.~K.}\ \bibnamefont
			{Seo}}\ and\ \bibinfo {author} {\bibfnamefont {E.~J.}\ \bibnamefont {Woo}},\
	}\enquote {\bibinfo {title} {Nonlinear inverse problems in imaging},}\ \
	(\bibinfo  {publisher} {John Wiley \& Sons, Ltd, Chichester, UK},\ \bibinfo
	{year} {2013})\ Chap.\ \bibinfo {chapter} {Electrical {I}mpedance
		{T}omography}\BibitemShut {NoStop}%
	\bibitem [{\citenamefont {Cultrera}\ \emph {et~al.}(2019)\citenamefont
		{Cultrera}, \citenamefont {Serazio}, \citenamefont {Zurutuza}, \citenamefont
		{Centeno}, \citenamefont {Txoperena}, \citenamefont {Etayo}, \citenamefont
		{Cordon}, \citenamefont {Redo-Sanchez}, \citenamefont {Arnedo}, \citenamefont
		{Ortolano},\ and\ \citenamefont {Callegaro}}]{cultrera2019mapping}%
	\BibitemOpen
	\bibfield  {author} {\bibinfo {author} {\bibfnamefont {A.}~\bibnamefont
			{Cultrera}}, \bibinfo {author} {\bibfnamefont {D.}~\bibnamefont {Serazio}},
		\bibinfo {author} {\bibfnamefont {A.}~\bibnamefont {Zurutuza}}, \bibinfo
		{author} {\bibfnamefont {A.}~\bibnamefont {Centeno}}, \bibinfo {author}
		{\bibfnamefont {O.}~\bibnamefont {Txoperena}}, \bibinfo {author}
		{\bibfnamefont {D.}~\bibnamefont {Etayo}}, \bibinfo {author} {\bibfnamefont
			{A.}~\bibnamefont {Cordon}}, \bibinfo {author} {\bibfnamefont
			{A.}~\bibnamefont {Redo-Sanchez}}, \bibinfo {author} {\bibfnamefont
			{I.}~\bibnamefont {Arnedo}}, \bibinfo {author} {\bibfnamefont
			{M.}~\bibnamefont {Ortolano}}, \ and\ \bibinfo {author} {\bibfnamefont
			{L.}~\bibnamefont {Callegaro}},\ }\href {\doibase 10.1038/s41598-019-46713-8}
	{\bibfield  {journal} {\bibinfo  {journal} {Sci. Rep.}\ }\textbf {\bibinfo
			{volume} {9}},\ \bibinfo {pages} {10655} (\bibinfo {year}
		{2019})}\BibitemShut {NoStop}%
	\bibitem [{\citenamefont {Cultrera}\ and\ \citenamefont
		{Callegaro}(2016)}]{cultrera2016electrical}%
	\BibitemOpen
	\bibfield  {author} {\bibinfo {author} {\bibfnamefont {A.}~\bibnamefont
			{Cultrera}}\ and\ \bibinfo {author} {\bibfnamefont {L.}~\bibnamefont
			{Callegaro}},\ }\href@noop {} {\bibfield  {journal} {\bibinfo  {journal}
			{IEEE Trans. Instrum. Meas.}\ }\textbf {\bibinfo {volume} {65}},\ \bibinfo
		{pages} {2101} (\bibinfo {year} {2016})}\BibitemShut {NoStop}%
	\bibitem [{\citenamefont {Adler}\ and\ \citenamefont
		{Lionheart}(2005)}]{adler2005eidors}%
	\BibitemOpen
	\bibfield  {author} {\bibinfo {author} {\bibfnamefont {A.}~\bibnamefont
			{Adler}}\ and\ \bibinfo {author} {\bibfnamefont {W.~R.}\ \bibnamefont
			{Lionheart}},\ }in\ \href@noop {} {\emph {\bibinfo {booktitle} {{6th
					Conference on Biomedical Applications of Electrical Impedance Tomography
					(London, UK, 22--24 June 2005)}}}}\ (\bibinfo {year} {2005})\BibitemShut
	{NoStop}%
\end{thebibliography}
%

\end{document}